\documentclass[11pt]{amsart}
\usepackage{amsmath,amsthm,amssymb}
\usepackage{enumerate}
\usepackage{graphicx}
\usepackage{cite}
\usepackage{comment}
\usepackage{oands}
\usepackage{tikz}
\usepackage{changepage}
\usepackage{bbm}
\usepackage{mathtools}
\usepackage[margin=1in]{geometry}
\usepackage{hyperref}
\usepackage[pagewise,mathlines]{lineno}
\usepackage{appendix}
\usepackage{multicol}
\usepackage{stmaryrd} 
\usepackage{enumitem}

\theoremstyle{plain}
\newtheorem{thm}{Theorem}

\hypersetup{
    colorlinks=false,
    linktocpage,
    }

\theoremstyle{definition}

\newtheorem{prob}[thm]{Problem}

\numberwithin{equation}{section}

\newcommand{\dsb}{\begin{adjustwidth}{2.5em}{0pt}
\begin{footnotesize}}
\newcommand{\dse}{\end{footnotesize}
\end{adjustwidth}}

\newcommand{\ssb}{\begin{adjustwidth}{2.5em}{0pt}}
\newcommand{\sse}{\end{adjustwidth}}

\newcommand{\aryb}{\begin{eqnarray*}}
\newcommand{\arye}{\end{eqnarray*}}
\def\alb#1\ale{\begin{align*}#1\end{align*}}
\def\allb#1\alle{\begin{align}#1\end{align}}
\newcommand{\eqb}{\begin{equation}}
\newcommand{\eqe}{\end{equation}}
\newcommand{\eqbn}{\begin{equation*}}
\newcommand{\eqen}{\end{equation*}}

\newcommand{\BB}{\mathbbm}

\newcommand{\op}{\operatorname}

\newcommand{\frk}{\mathfrak}

\newcommand{\ep}{\varepsilon}
\newcommand{\rta}{\rightarrow}

\newcommand{\wt}{\widetilde}
 
\newcommand{\mcl}{\mathcal}

\let\originalleft\left
\let\originalright\right
\renewcommand{\left}{\mathopen{}\mathclose\bgroup\originalleft}
\renewcommand{\right}{\aftergroup\egroup\originalright}

\title{Random surfaces and Liouville quantum gravity} 
\author{Ewain Gwynne}

\begin{document}

\begin{abstract}
Liouville quantum gravity (LQG) surfaces are a family of random fractal surfaces which can be thought of as the canonical models of random two-dimensional Riemannian manifolds, in the same sense that Brownian motion is the canonical model of a random path. LQG surfaces are the continuum limits of discrete random surfaces called random planar maps. 
In this expository article, we discuss the definition of random planar maps and LQG, the sense in which random planar maps converge to LQG, and the motivations for studying these objects. We also mention several open problems.
We do not assume any background knowledge beyond that of a second-year mathematics graduate student. 
\end{abstract}

\maketitle

What is the most natural way of choosing a random surface (two-dimensional Riemannian manifold)?
If we are given a finite set $X$, the easiest way to choose a random element of $X$ is uniformly, i.e., by assigning equal probability to each element of $X$.  More generally, if we are given a set $X\subset\BB R^n$ with finite, positive Lebesgue measure, the simplest way of choosing a random element of $X$ is by sampling from Lebesgue measure normalized to have total mass one. 
However, the space of all surfaces is infinite dimensional for any reasonable notion of dimension, so it is not immediately obvious whether there is a canonical way of choosing a random surface. 

Nevertheless, there is a class of canonical models of ``random surfaces" called \emph{Liouville quantum gravity (LQG) surfaces}.  
The reason for the quotations is that LQG surfaces are not Riemannian manifolds in the literal sense since they are too singular to admit a smooth structure.
Instead, LQG surfaces are defined as random topological surfaces equipped with a measure, a metric, and a conformal structure.
These surfaces are fractal in the sense that the Hausdorff dimension of an LQG surface, viewed as a metric space, is strictly bigger than 2. 

LQG surfaces have a rich geometric structure which is still not fully understood (see Section~\ref{sec-open-problems}). 
Furthermore, such surfaces are important in statistical mechanics, string theory, and conformal field theory and have deep connections to other mathematical objects such as Schramm-Loewner evolution~\cite{schramm0} (see, e.g.,~\cite{shef-zipper}), random matrix theory (see, e.g.,~\cite{webb-gmc-rmt}), and random planar maps.

\section{Discrete random surfaces}
\label{sec-rpm}

Before we discuss random surfaces, let us first, by way of analogy, consider the simpler problem of finding a canonical way to choose a random curve in the plane.
As in the case of surfaces, the space of all planar curves is infinite-dimensional and does not admit a canonical probability measure in an obvious way.
To get around this, we discretize the problem. 
Let us consider for each $n\in\BB N$ the set of all nearest-neighbor paths in the integer lattice $\BB Z^2$ with $n$ steps. 
This is a finite set (of cardinality $4^n$), so we can choose a uniformly random element $S_n$ from it.
This discrete random path $S_n$ is called the \emph{simple random walk}. 

By linearly interpolating, we may view $S_n$ as a curve from $[0,n]$ to $\BB R^2$. 
There is a classical theorem in probability due to Donsker which states that the re-scaled random walk paths $t\mapsto n^{-1/2} S_n(n t)$ converge in distribution (with respect to the uniform topology) to a limiting random continuous curve called \emph{Brownian motion}, which can be thought of as the canonical random planar curve. 
Brownian motion has a much richer structure than an ordinary smooth curve. A Brownian motion curve is nowhere differentiable and crosses itself in every interval of time. It has zero Lebesgue measure, but Hausdorff dimension 2. 

\begin{figure} 
\begin{center}
\includegraphics[scale=.6]{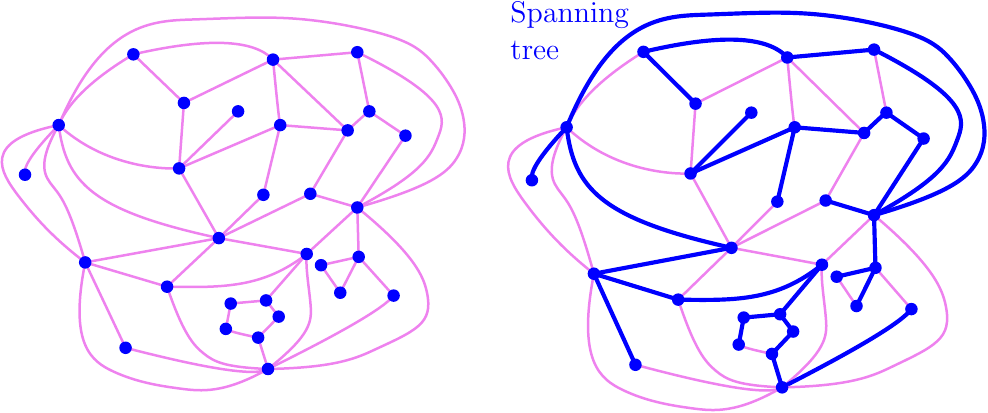} 
\caption{\label{fig-planar-map} 
\textbf{Left:} A planar map. \textbf{Right:} A planar map decorated by a spanning tree.
}
\end{center}
\end{figure}

One can define a canonical random surface via a similar approach. Let us first define the discrete random surfaces which we will consider. 
A \emph{planar map} is a graph (multiple edges and self-loops allowed) embedded into the plane $\BB C$ in such a way that no two edges cross, viewed modulo orientation-preserving homeomorphisms $\BB C\rta\BB C$. 
Planar maps have vertices, edges, and faces, but these objects do not correspond to particular points / sets in $\BB C$ since we do not specify a particular embedding. 
See Figure~\ref{fig-planar-map} for an illustration of a planar map.
 
We can think of a planar map as a discrete random surface, where each face is equipped with the flat Riemannian metric coming from a polygon in $\BB C$ with the appropriate number of edges and unit side length. 
Equivalently, a planar map is obtained by starting with a collection of polygons with unit side length and identifying pairs of their sides (in such a way that the identification between any two edges is a Euclidean isometry) to produce a surface, subject to the constraint that there are no holes or handles. 
The surface obtained in this way always has the topology of the sphere, but we can similarly obtain random surfaces with other topologies (we are primarily interested in the local geometry of random surfaces rather than their topology). 

For each $n\in\BB N$, there are only finitely many planar maps with $n$ edges. Hence it makes sense to choose such a map uniformly at random.
One can also consider uniform planar maps with local constraints, such as triangulations (resp.\ quadrangulations), which are required to have three (resp.\ four) edges on the boundary of each face. 

It is also natural to consider planar maps weighted by some sort of additional structure on the map.
Indeed, suppose, for example, that we are interested in planar maps $M$ decorated by a \emph{spanning tree}, i.e., a connected subgraph of $M$ which contains every vertex of $M$ and has no cycles.
Then we might want sample a uniform random pair $(M,T)$ consisting of a planar map with $n$ edges and a spanning tree on it.\footnote{There is some ambiguity when counting pairs $(M,T)$ since there can be automorphisms $M \rta M$ which do not fix $T$. In practice, this ambiguity is removed by specifying a distinguished edge $e$ of $M$: there are no non-trivial automorphisms of $M$ which fix $e$.}
In this case, the marginal distribution of $M$ is not uniform: rather, the probability that $M$ is equal to any fixed planar map $\frk M$ with $n$ edges is proportional to the number of possible spanning trees of $\frk M$. 
In a similar vein, one might want to look at planar maps sampled with probability proportional to the number of certain types of orientations on the edges of $M$, or the partition function of a statistical mechanics model on $M$ (like the Ising model or the Fortuin-Kasteleyn model).

Since planar maps can be thought of as discrete surfaces, it is natural to expect that random planar maps converge, in some sense, to limiting random surfaces as the total number of edges tends to $\infty$. In other words, if we sample a large random planar map and ``zoom out" so that we see only its large-scale structure, we should get something which looks like some sort of (continuum) random surface. 
For a certain class of random planar maps, including uniform planar maps and planar maps sampled with probability proportional to the partition function of a statistical mechanics model with parameters set to the critical values, the particular types of surfaces which arise in this way are LQG surfaces, as discussed above and defined precisely in Section~\ref{sec-lqg}. 
We will discuss the precise sense in which this convergence occurs in Section~\ref{sec-conv}.

\section{Liouville quantum gravity} 
\label{sec-lqg}

\subsection{Isothermal coordinates}
\label{sec-isothermal}
 
To define LQG surfaces, we first recall some facts from the theory of deterministic surfaces. 
Suppose $\mcl S$ is a continuously differentiable surface, i.e., in local coordinates $\mcl S$ can be represented by $E\,dx^2 + F\,dx\, dy + G \,dy^2$ for some continuously differentiable functions $E,F,G$ of the parameters $  (x,y)$. 
We will be primarily interested in the local geometry of $\mcl S$, so by possibly replacing $\mcl S$ by an open subset of $\mcl S$, we can assume that $\mcl S$ is homeomorphic to the open unit disk $\BB D\subset\BB R^2$. 
A standard theorem in Riemannian geometry (see, e.g.,~\cite{chern-isothermal-coord}) 
asserts that, at least locally, $\mcl S$ can be parametrized by \emph{isothermal coordinates}. This means that we can parametrize $\mcl S$ by coordinates $(x,y)$ in $\BB D$ in such a way that the Riemannian metric tensor takes the form $e^{h(z)} \, (dx^2 + dy^2)$ for some continuously differentiable function $h : \BB D\rta \BB R$.
Here, $z = x+i y$ and $dx^2+dy^2$ is the Euclidean Riemannian metric tensor.  
 
It is easy to describe areas and distances with respect to isothermal coordinates. For a Lebesgue measurable set $A\subset\BB D$, the area of the corresponding subset of the surface $\mcl S$ is given by 
\eqb \label{eqn-surface-area}
\int_A e^{h(z)} \,d^2 z  ,
\eqe  
where $d^2 z = dx\,dy$ is two-dimensional Lebesgue measure.
The $\mcl S$-distance between any two points in $z,w \in \BB D$ is given by
\eqb \label{eqn-surface-dist}
\inf_{P : z\rta w} \int_a^b e^{h(P(t))/2} |P'(t)| \,dt ,
\eqe
where the inf is over all piecewise continuously differentiable paths $P: [a,b] \rta \BB D$ from $z$ to $w$.

\subsection{The Gaussian free field}
\label{sec-gff}

We want to define an LQG surface as a random surface parametrized by isothermal coordinates by making a random choice of $h$. 
Since LQG surfaces should describe the large-scale behavior of random planar maps, by analogy with the central limit theorem a natural first guess is that $h$ should be a ``standard Gaussian random variable" taking values in the space of differentiable functions on $\BB D$. To explain what this means, suppose that we are given a finite-dimensional Hilbert space $\mcl H$ and let $\{x_1,\dots,x_n\}$ be an orthonormal basis for $\mcl H$.
We can define a \emph{standard Gaussian random variable} on $\mcl H$ by 
\eqb \label{eqn-standard-gaussian}
  \sum_{j=1}^n \alpha_j x_j 
\eqe
where the $\alpha_j$'s are i.i.d.\ standard Gaussian random variables on $\BB R$, i.e., they are sampled from the probability measure with density $\frac{1}{\sqrt{2\pi}} e^{-x^2/2}$. 

For an open domain $U\subset\BB C$, consider the infinite-dimensional Hilbert space $\mcl H(U)$ which is the Hilbert space completion of the space of smooth, compactly supported functions on $\BB D$ with respect to the \emph{Dirichlet inner product}
\eqb
(f,g)_\nabla = \frac{1}{2\pi} \int_{U} \nabla f(z) \cdot \nabla g(z) \,d^2 z ,
\eqe
where $\nabla$ denotes the gradient and $\cdot$ denotes the dot product.
This is sometimes called the \emph{first order Sobolev space} on $U$, with zero boundary conditions. 

The \emph{Gaussian free field (GFF)} is the standard Gaussian random variable on $\mcl H(U)$. 
That is, let $\{f_j\}_{j\in\BB N}$ be an orthonormal basis for $\mcl H(U)$. 
By analogy with~\eqref{eqn-standard-gaussian}, we define the GFF $h$ by
\eqb \label{eqn-gff-def}
h := \sum_{j=1}^\infty \alpha_j f_j 
\eqe
where the $\alpha_j$'s are i.i.d.\ standard Gaussian random variables. 
The sum~\eqref{eqn-gff-def} a.s.\ does not converge pointwise, so the GFF does not have well-defined pointwise values.

However, it is not hard to show that the GFF makes sense as a random distribution (generalized function) on $\mcl H(U)$~\cite{shef-gff}.
This means that for any $f\in\mcl H(U)$ the Dirichlet inner product ``$(h,f)_\nabla$" and the $L^2$ inner product $(h,\phi) = ``\int_U h(z) f(z) \,d^2 z"$ are well-defined.  
The random variables $(h,f)_\nabla$ and $(h,g)_\nabla$ for $f,g\in\mcl H(U)$ are jointly centered Gaussian with covariances $\op{Cov}((h,f)_\nabla , (h,g)_\nabla) = (f,g)_\nabla$.

Although we will primarily be interested in the two-dimensional case, we remark that the above construction of the GFF also makes sense in other dimensions. 
In dimension 1, one gets a Brownian bridge (a one-dimensional Brownian motion defined on an interval and conditioned to be zero at the endpoints). Hence the GFF can be seen as a generalization of Brownian motion with two time dimensions. 
In dimension at least three, the GFF is in some sense rougher, i.e., further from being a function, than in dimension 2 and hence regularization procedures such as the ones described below do not give interesting objects. Rather, it is more natural to consider so-called \emph{log-correlated} Gaussian fields~\cite{lgf-survey} (the GFF is log-correlated only in dimension 2).

 
\subsection{Liouville quantum gravity surfaces} 
\label{sec-lqg-def}

Let $h$ be the GFF on a domain $U\subset\BB C$ and let $\gamma\in (0,2)$.  
Recalling the above discussion on isothermal coordinates, we define  the \emph{$\gamma$-Liouville quantum gravity (LQG) surface} associated with $(U,h)$ to be the random surface parametrized by $U$ with Riemannian metric tensor ``$e^{\gamma h} \,(dx^2 + dy^2)$". 
This definition does not make literal sense since $h$ is not a function, so it cannot be exponentiated. 
However, one can define LQG surfaces rigorously via regularization procedures which we will describe shortly.

The parameter $\gamma$ controls the ``roughness" of the surface: $\gamma=0$ corresponds to a smooth surface, and larger values of $\gamma$ correspond to surfaces which are more fractal and less Euclidean-like. 
The parameter $\gamma$ is also related to the type of random planar map model under consideration. 
The case when $\gamma=\sqrt{8/3}$ (sometimes called ``pure gravity") describes the large-scale behavior of uniform planar maps. Other values of $\gamma$ (``gravity coupled to matter") correspond to random planar maps sampled  with probability proportional to some sort of additional structure on the map, such as the number of spanning trees ($\gamma = \sqrt{2}$) or the Ising model partition function ($\gamma = \sqrt{3}$).

LQG surfaces were first defined non-rigorously in the physics literature by Polyakov~\cite{polyakov-qg1} in the context of bosonic string theory.\footnote{Roughly speaking, for $\mathbf c \in\BB N$ an evolving string in $\BB R^{\mathbf c-1}$ traces out a two-dimensional surface embedded in space-time $\BB R^{\mathbf c-1} \times \BB R$, called a \emph{world sheet}. Polyakov wanted to develop a theory of integrals over all possible surfaces embedded in $\BB R^{\mathbf c}$ as a string-theoretic generalization of the Feynman path integral (which is an integral over all possible paths). To do this one needs to define a probability measure on surfaces.
The most natural way of doing this turns out to only work when the ``dimension of the space into which the surface is embedded" (a.k.a.\ the central charge) $\mathbf c$ lies in $ (-\infty,1]$, in which case the desired probability measure is $\gamma$-LQG with $\mathbf c = 25-6(2/\gamma+\gamma/2)^2$. }
We refer to~\cite{shef-kpz} for an extensive list of references to the physics literature on LQG. 
Certain special LQG surfaces are the topic of study in Liouville conformal field theory, the simplest non-rational conformal field theory.
See, e.g.,~\cite{dkrv-lqg-sphere,grv-higher-genus,krv-dozz} and the references therein for rigorous works on LQG from the conformal field theory perspective. 

We now explain how to make rigorous sense of LQG surfaces as random metric measure spaces. 
The basic idea is to consider a family of continuous functions $\{h_\ep\}_{\ep > 0}$ which approximate the GFF as $\ep\rta 0$, define approximate notions of area and distance by replacing $h$ with a multiple of $h_\ep$ in~\eqref{eqn-surface-area} and~\eqref{eqn-surface-dist}, then send $\ep \rta 0$ and re-normalize appropriately to get a measure and metric associated with $\gamma$-LQG. 
One possible choice of $h_\ep$ is the convolution with the heat kernel, 
\eqb
h_\ep(z) =  \int_U h(z) p_\ep(z,w) \,d^2w \quad \text{for} \quad p_\ep(z,w) := \frac{1}{\pi \ep^2} e^{-|z-w|^2/\ep^2} 
\eqe
where the integral is interpreted as the distributional pairing of $h$ with $p_\ep(z,\cdot)$. 
Note that $p_\ep(z,w)\,d^2w$ approximates a point mass at $z$ as $\ep\rta 0$, so $h_\ep$ is close to $h$ (e.g., in the distributional sense) when $\ep$ is small. 
Other possible choices for $h_\ep$ include convolutions with other mollifiers, averages over circles, truncated versions of the orthonormal basis expansion~\eqref{eqn-gff-def}, etc.

\begin{figure}[ht!]
\begin{center}
\includegraphics[scale=.6]{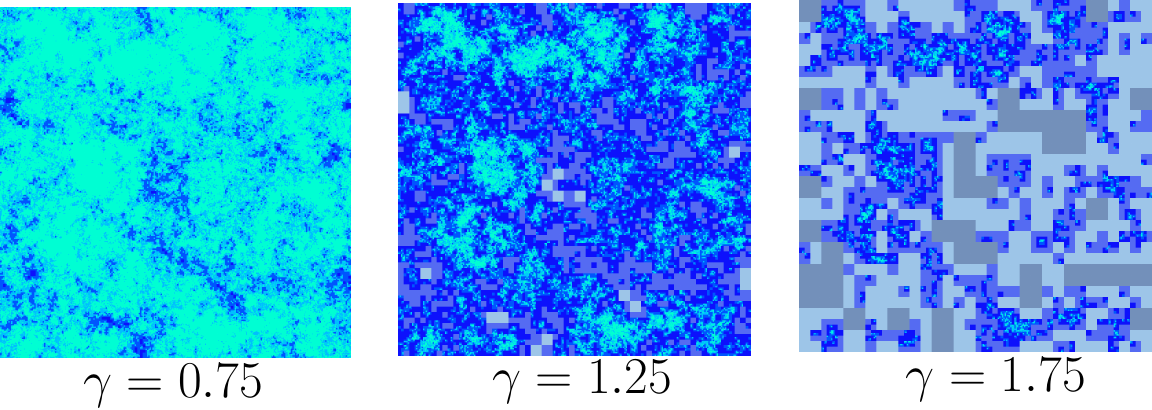} 
\caption{\label{fig-measure-sim}
Simulations of the $\gamma$-LQG measure on the unit square produced by J.\ Miller.
The square is divided into dyadic sub-squares which all have approximately the same $\gamma$-LQG mass.
Squares are colored according to their Euclidean size.
Note that as $\gamma$ increases, the Euclidean sizes of these squares become more variable. 
}
\end{center}
\end{figure}
 
\subsection{The $\gamma$-LQG area measure}

For several possible choices of $\{h_\ep\}_{\ep>0}$, one can define the \emph{$\gamma$-LQG area measure} as the a.s.\ limit
\eqb \label{eqn-lqg-measure}
\mu_h = \lim_{\ep\rta 0} \ep^{\gamma^2/2} e^{\gamma h_\ep(z)} \,d^2 z ,
\eqe
with respect to the weak topology on measures on $U$, where $d^2 z$ denotes Lebesgue measure as above. 
Indeed, this construction is a special case of a more general theory of regularized random measures called \emph{Gaussian multiplicative chaos}, which was initiated in the work of Kahane~\cite{kahane}; see~\cite{rhodes-vargas-review,berestycki-gmt-elementary,aru-gmc-survey} for surveys of this theory. 
Several important properties of the measure (including the convergence of the circle average approximation and the so-called \emph{KPZ formula}) were estalbished by Duplantier and Sheffield~\cite{shef-kpz}. 
The $\gamma$-LQG measure has no point masses and assigns positive mass to every open subset of $U$, but it is mutually singular with respect to the Lebesgue measure. In fact, it assigns full mass to a set of Hausdorff dimension $2-\gamma^2/2$ (see, e.g., ~\cite[Section 3.3]{shef-kpz} and~\cite{hmp-thick-pts}). 
See Figure~\ref{fig-measure-sim} for simulations of the $\gamma$-LQG measure.

In the case when $\gamma \geq 2$, the limit~\eqref{eqn-lqg-measure} is identically zero, which is why we restrict to $\gamma \in (0,2)$ when defining an LQG surface. 
For $\gamma=2$ one can construct a measure with similar properties but an additional logarithmic correction is needed in the scaling factor; see~\cite{shef-deriv-mart,shef-renormalization}. For $\gamma > 2$, one can make sense of the $\gamma$-LQG measure as a purely atomic measure (i.e., a countable sum of point masses) which is closely related to $\gamma'$-LQG for $\gamma'  = 4/\gamma \in (0,2)$~\cite{dup-dual-lqg,bjrv-gmt-duality,rhodes-vargas-review,wedges}. 

\subsection{The $\gamma$-LQG metric}

The $\gamma$-LQG metric $D_h$ can be constructed in an analogous way to~\eqref{eqn-lqg-measure}, but the proof that the approximations converge is much more involved than in the case of the measure. Intuitively, the reason for this is that if we replace $h$ by a multiple of $ h_\ep$ in~\eqref{eqn-surface-dist}, then the near-minimal paths could in principle be very different for different values of $\ep$ (although one gets \emph{a posteriori} that this is not the case). 

Before describing the construction of the metric we need to introduce an exponent $d_\gamma > 2$ which plays a fundamental role in the study of LQG distances and which can be defined in several equivalent ways. 
For example, it was shown in~\cite{dg-lqg-dim}, building on~\cite{dzz-heat-kernel,ghs-map-dist} that for certain random planar maps in the $\gamma$-universality class, a graph distance ball of radius $r \in \BB N$ in the map typically has of order $r^{d_\gamma}$ vertices.
Once the $\gamma$-LQG metric $D_h$ is constructed, it is possible to show that $d_\gamma$ is the Hausdorff dimension of the metric space $(U,D_h)$~\cite{gp-kpz}.

It can be shown using special symmetries for uniform planar maps or $\sqrt{8/3}$-LQG that $d_{\sqrt{8/3}}=4$.
However, $d_\gamma$ is not known (even at a heuristic level) for $\gamma \in (0,2) \setminus \{\sqrt{8/3}\}$; determining its value is one of the most important open problems in the theory of LQG. See Section~\ref{sec-open-problems} for more on $d_\gamma$. 
 
To approximate LQG distances, we let $D_h^\ep(z,w)$ for $z,w\in U$ and $\ep > 0$ be the metric defined as in~\eqref{eqn-surface-dist} with $  (2\gamma/d_\gamma) h_\ep$ in place of $h$, i.e.,
\eqb \label{eqn-lfpp-def}
D_h^\ep(z,w) = \inf_{P : z\rta w} \int_a^b e^{\frac{\gamma}{d_\gamma} h_\ep(P(t)) } |P'(t)| \,dt .
\eqe
The reason why we have $\gamma/d_\gamma$ instead of $\gamma/2$ in~\eqref{eqn-lfpp-def} is as follows. 
Since $d_\gamma$ is the dimension of the $\gamma$-LQG surface, scaling $\gamma$-LQG areas by $C> 0$ should correspond to scaling $\gamma$-LQG distances by $C^{1/d_\gamma}$. 
By~\eqref{eqn-lqg-measure}, scaling areas by $C$ corresponds to adding the constant $\frac{1}{\gamma} \log C$ to $h$.
By~\eqref{eqn-lfpp-def}, this scales $D_h^\ep$ by $C^{1/d_\gamma}$, as desired. 
 
It was shown by Ding-Dub\'edat-Dunlap-Falconet~\cite{dddf-lfpp} that there are constants $\{\frk a_\ep\}_{\ep > 0}$ such that the re-scaled metrics $\frk a_\ep^{-1} D_h^\ep$ are tight with respect to the local uniform topology on $U$, and every subsequential limit is bi-H\"older continuous with respect to the Euclidean metric on $U$.
Building on this and~\cite{local-metrics,lqg-metric-estimates,gm-confluence}, Gwynne and Miller~\cite{gm-uniqueness} showed that in fact $\frk a_\ep^{-1} D_h^\ep$ converges in probability (not just subsequentially) to a limiting metric $D_h$ which is defined to be the $\gamma$-LQG metric. 
They also proved an axiomatic characterization of $D_h$ which implies that it is the only possible metric associated with $\gamma$-LQG.

\begin{figure}[ht!]
\begin{center}
\includegraphics[scale=.8]{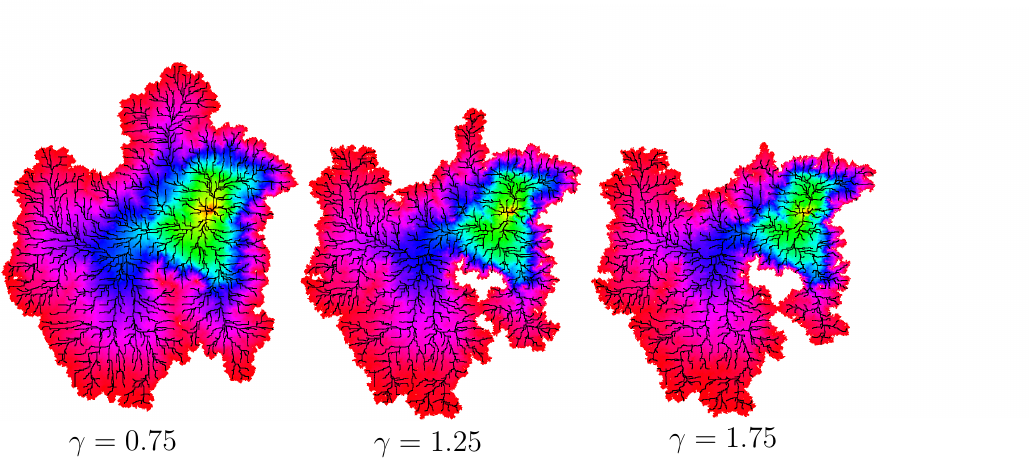} 
\caption{\label{fig-metric-sim}
Simulations of $\gamma$-LQG metric balls w.r.t.\ the same GFF instance, produced by J.\ Miller.
The colors indicate distances to the center of the ball.
Geodesics from points in a grid back to the center point are shown in black.
Note that these geodesics have a tree-like structure: unlike geodesics for a smooth Riemannian metric, LQG geodesics with different starting points and targeted at 0 merge into one another before reaching 0.
This was proven to be the case in~\cite{gm-confluence}. 
}
\end{center}
\end{figure}

The metric $D_h$ induces the same topology on $U$ as the Euclidean metric, but the Hausdorff dimension of the metric space $(U,D_h)$ is $d_\gamma > 2$. Moreover, many of its geometric properties (e.g., scaling properties and the behavior of geodesics) are quite different from those of the Euclidean metric or indeed any smooth Riemannian metric on $U$. See Figure~\ref{fig-metric-sim} for simulations of LQG metric balls. 

There is also an earlier construction of the LQG metric in the special case when $\gamma=\sqrt{8/3}$ due to Miller and Sheffield~\cite{lqg-tbm1,lqg-tbm2,lqg-tbm3}. 
This construction does not use the direct regularization procedure~\eqref{eqn-lfpp-def} and instead is based on special symmetries for $\sqrt{8/3}$-LQG.
The Miller-Sheffield construction does not work for $\gamma\not=\sqrt{8/3}$, but it gives additional information about $\sqrt{8/3}$-LQG metric (e.g., its connection to random planar maps and certain Markov properties of metric balls) which is not apparent from the construction of~\cite{gm-uniqueness}.  
The two constructions give the same metric by the characterization theorem from~\cite{gm-uniqueness}.

\subsection{Conformal coordinate change}

Just like for deterministic surfaces, it is possible to parametrize LQG surfaces in different ways. 
Suppose $\phi : \wt U \rta U$ is a conformal (i.e., bijective and holomorphic) map.
Let $h$ be the GFF on $U$ and let 
\eqb \label{eqn-lqg-coord}
\wt h := h\circ \phi + Q\log |\phi'| \quad \text{where} \quad Q := \frac{2}{\gamma} + \frac{\gamma}{2} .
\eqe
Then $\wt h$ is a random distribution on $\wt U$. It is shown in~\cite{shef-kpz,gm-coord-change} that the $\gamma$-LQG area measures and metrics associated with $h$ and $\wt h$ are a.s.\ related by $\mu_h(\phi(A)) = \mu_{\wt h}(A)$, for each Borel measurable set $A\subset \wt U$ and $D_h(\phi(z),\phi(w)) = D_{\wt h}(z,w)$ for each $z,w\in \wt U$.
We think of $(U,h)$ and $(\wt U , \wt h)$ as representing different parametrizations of the same LQG surface.
The coordinate change relation for $\mu_h$ and $D_h$ shows that these objects depend only on the LQG surface, not on the choice of parametrization.

\section{LQG as the limit of random planar maps}
\label{sec-conv}

We now discuss the senses in which random planar maps should converge to LQG surfaces, and the extent to which each type of convergence has been proven.

\begin{figure}[ht!]
\begin{center}
\includegraphics[scale=.35]{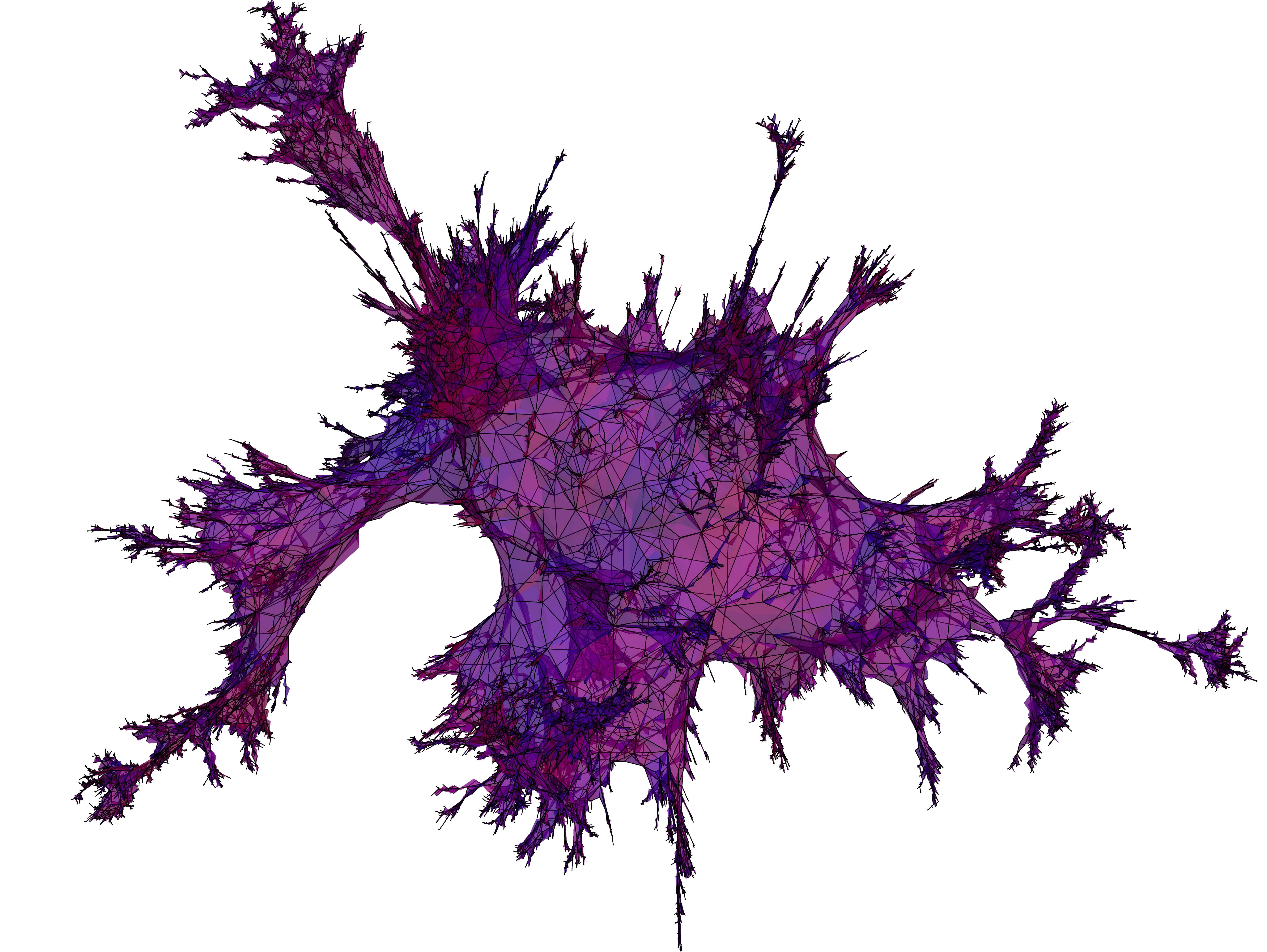}\hspace{10pt} 
\includegraphics[scale=.1]{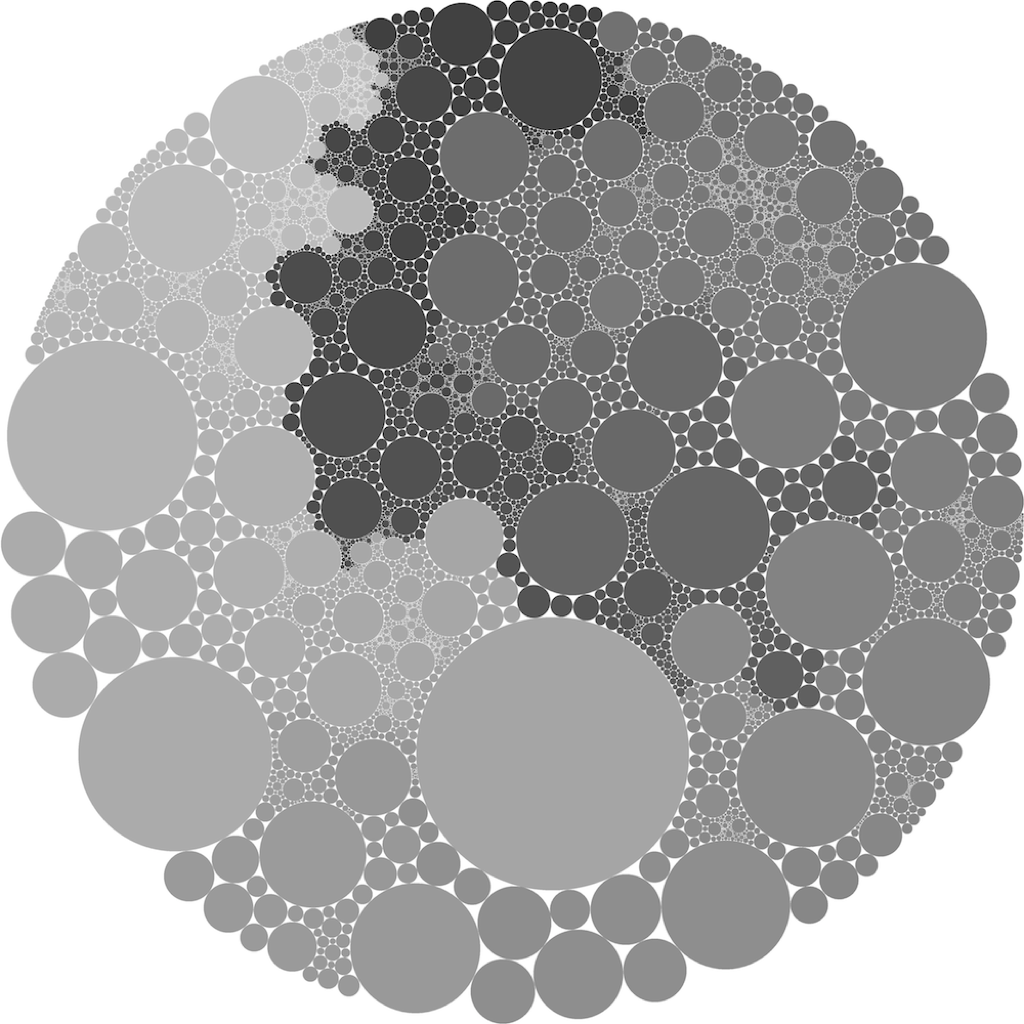}   \hspace{10pt} 
\includegraphics[scale=.1]{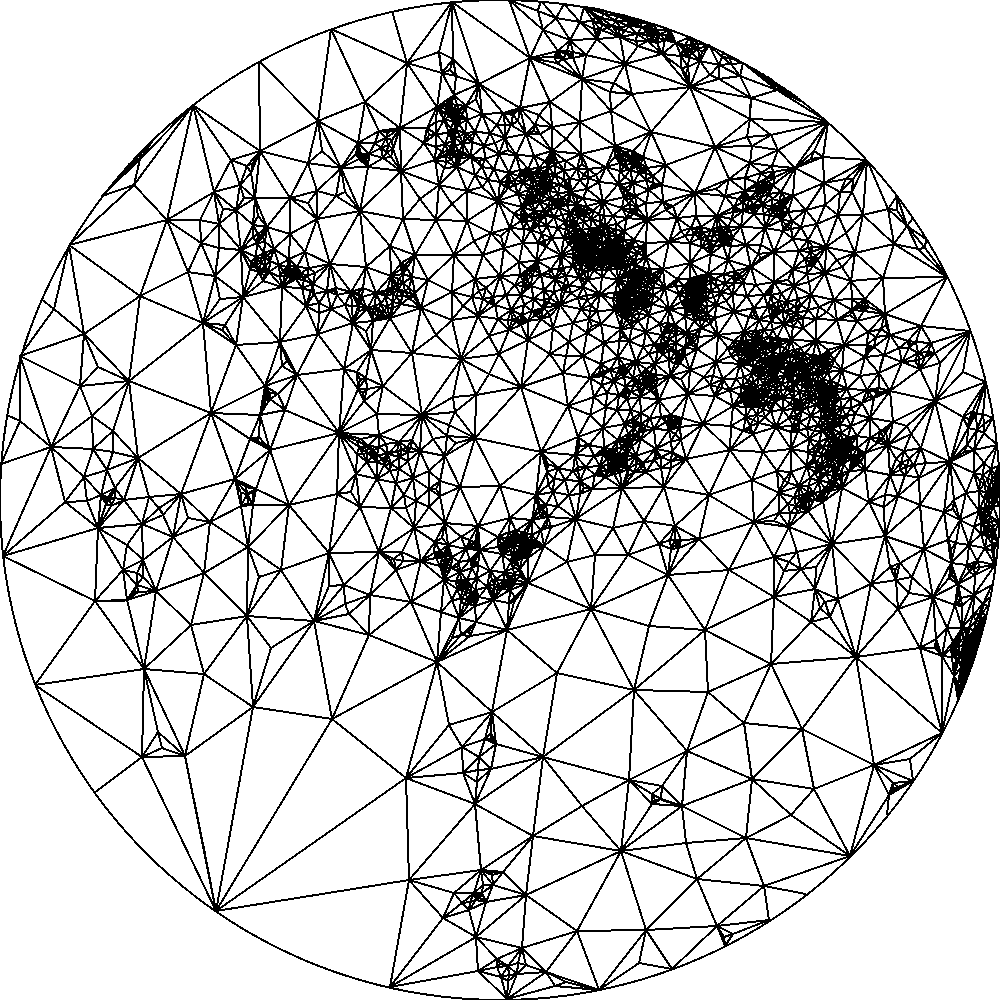}  
\caption{\label{fig-map-embedding} 
\textbf{Left:} Simulation of a large uniform quadrangulation embedded into $\BB R^3$ in such a way that the embedding is in some sense as close as possible to being an isometry, made by J.\ Bettinelli. The Gromov-Hausdorff limit of these triangulations is a $\sqrt{8/3}$-LQG surface with the topology of the sphere.
\textbf{Middle:}  Simulation of a circle packing of a uniform triangulation made by J.\ Miller. 
\textbf{Right:} Simulation of the Tutte embedding of an instance of the $\sqrt 2$-mated CRT map, made by J.\ Miller. The counting measure on vertices approximates the $\sqrt{2}$-LQG measure when the number of vertices is large.
}
\end{center}
\end{figure}

\subsection{Gromov-Hausdorff convergence} 
\label{sec-gh}

Let $\mcl K$ be the set of all compact metric spaces. The \emph{Gromov-Hausdorff (GH) distance} on $\mcl K$ is the metric on $\mcl K$ defined by
\eqb
\BB d^{\op{GH}}\left((X_1,d_1) , (X_2,d_2) \right)
:= \inf_{(Y,D) , \iota_1,\iota_2} \BB d^{\op{H}}\left(\iota_1(X_1) , \iota_2(X_2) \right)
\eqe
where $\BB d^{\op{H}}$ denotes the Hausdorff distance on compact subsets of $Y$ and the infimum is over all compact metric spaces $(Y,D)$ and isometric embeddings $\iota_1: (X_1,d_1) \rta (Y,D)$ and $\iota_2 : (X_2,d_2) \rta (Y,D)$.
A planar map can be viewed as a compact metric space equipped with its graph distance (see Figure~\ref{fig-map-embedding}, left).
One can then ask whether large random planar maps, with their graph distance re-scaled appropriately, converge in distribution to $\gamma$-LQG surfaces w.r.t.\ the GH topology.
We can additionally equip the map with its re-scaled counting measure and ask for convergence in the \emph{Gromov-Hausdorff-Prokhorov (GHP) topology}, the analog of the GH topology for metric measure spaces.

So far, GH or GHP convergence has only been established for uniform planar maps (including uniform maps with local constraints), which we recall correspond to $\gamma=\sqrt{8/3}$. 
The first such convergence results were obtained independently by Le Gall~\cite{legall-uniqueness} and Miermont~\cite{miermont-brownian-map}.
They showed that certain types of uniform planar maps (with graph distances re-scaled by $n^{-1/4}$ and the counting measure on vertices scaled by $1/n$) converge in the GHP sense to the \emph{Brownian map}, a random metric measure space which can be constructed from a continuum random tree via an explicit metric quotient procedure.
Miller and Sheffield~\cite{lqg-tbm1,lqg-tbm2,lqg-tbm3} showed that a certain special $\sqrt{8/3}$-LQG surface called the \emph{quantum sphere}, has the same distribution as the Brownian map viewed as a metric measure space modulo measure-preserving isometries. 
Hence uniform planar maps converge to $\sqrt{8/3}$-LQG surfaces in the metric space sense. However, the proofs in~\cite{legall-uniqueness,miermont-brownian-map,lqg-tbm1,lqg-tbm2,lqg-tbm3} rely on special symmetries and combinatorial miracles~\cite{schaeffer-bijection,bdg-bijection,angel-peeling} for uniform planar maps and $\sqrt{8/3}$-LQG (these are the same sorts of miracles which are needed to show $d_{\sqrt{8/3}}=4$). So, the convergence proofs do not extend to non-uniform maps and $\gamma\not=\sqrt{8/3}$. 

\subsection{Embedding convergence} 
\label{sec-embedding}

Although a planar map is defined only modulo orientation-preserving homeomorphisms of $\BB C$, there are various ways of \emph{embedding} the map into $\BB C$, i.e., associating each vertex (resp.\ edge) with a point (resp.\ curve) in $\BB C$ in such a way that no two edges cross.
Examples include \emph{circle packing}, where the planar map is realized as the tangency graph of a collection of circles in the plane (Figure~\ref{fig-map-embedding}, middle); and \emph{Tutte embedding} (a.k.a.\ harmonic embedding or baycentric embedding), which is defined by the condition that the position of each vertex is the average (barycenter) of the positions of its neighbors (Figure~\ref{fig-map-embedding}, right). 
Once we have chosen an embedding, we can ask, e.g., whether the counting measure on the vertices of the planar map, re-scaled by the total number of vertices, converges in distribution (w.r.t.\ the weak topology on $\BB C$) to a variant of the $\gamma$-LQG area measure. 

So far, there are three results establishing this type of convergence for random planar maps.
The first result~\cite{gms-tutte} establishes embedding convergence for a one-parameter family of random planar maps called \emph{mated-CRT maps}, one for each $\gamma \in (0,2)$, embedded via the Tutte embedding. 
Mated-CRT maps are constructed by gluing together pairs of correlated continuum random trees (CRT's), and are more directly connected to LQG than other types of random planar maps due to the results of~\cite{wedges}. We will not define these maps precisely here. 

Embedding convergence toward $\sqrt{8/3}$-LQG has also been established for uniform triangulations embedded via the so-called \emph{Cardy embedding}, which is defined using crossing probabilities for percolation on the map~\cite{hs-cardy-embedding}; and for the Poisson-Voronoi approximation of the Brownian map under the Tutte embedding~\cite{gms-poisson-voronoi}. In the case of uniform triangulations, it is in fact shown in~\cite{hs-cardy-embedding} (see also~\cite{ghs-metric-peano}) that one has convergence in the GHP sense (as discussed in Section~\ref{sec-gh}) and the mating-of-trees sense (as discussed in Section~\ref{sec-peano} below) \emph{simultaneously} with the convergence of the embedded map.
It is a major open problem to establish these types of convergence results for other random planar map models and/or other embeddings (see Problem~\ref{prob-conv}).

\subsection{Mating-of-trees convergence}
\label{sec-peano}

There are several combinatorial bijections which encode a random planar map decorated by some additional structure by means of a random walk on $\BB Z^2$.
The simplest example of such a bijection is the Mullin bijection~\cite{mullin-maps,bernardi-maps}, which encodes a planar map decorated by a spanning tree by a nearest-neighbor random walk on $\BB Z^2$.
There are other such bijections, with different walks, that encode planar maps decorated by e.g., percolation~\cite{bernardi-dfs-bijection,bhs-site-perc}, bipolar orientations~\cite{kmsw-bipolar}, or the Fortuin-Kasteleyn model~\cite{shef-burger}.
These bijections are called \emph{mating-of-trees} bijections since the planar map is constructed from the walk by gluing together, or mating, the discrete random trees associated with the two coordinates of the walk. In the case of the Mullin bijection, the two trees are the spanning tree on the map and its corresponding dual spanning tree.

For random planar maps which are related to $\gamma$-LQG, it can be shown that the encoding walk converges in distribution to a two-dimensional Brownian motion $Z = (L,R)$ with the correlation between the coordinates $L_t$ and $R_t$ given by $-\cos(\pi\gamma^2/4)$ for each time $t$. 
A fundamental theorem of Duplantier, Miller, and Sheffield~\cite{wedges} shows that one can construct a $\gamma$-LQG surface from this two-dimensional Brownian motion $Z$ via a continuum analog of a mating-of-trees bijection. 
Hence the convergence of the encoding walks in the discrete mating-of-trees bijections toward $Z$ can be viewed as a convergence statement for random planar maps toward LQG surfaces in a certain topology: the one where two surfaces are close if their encoding paths are close. 
This type of convergence is referred to as \emph{mating-of-trees} or \emph{peanosphere} convergence.

Mating-of-trees convergence is not strictly weaker than Gromov-Hausdorff convergence or embedding convergence, but it is arguably less natural than these other modes of convergence.  
However, mating-of-trees convergence can still be used to extract a substantial amount of useful information about the random planar map. This includes scaling limits for various functionals of the map and the computation of exponents related to partition functions, graph distances, random curves on the map, etc. Moreover, for many types of random planar maps mating-of-trees convergence is the only scaling limit result available, and it can sometimes be used as an intermediate step in proving one of the other types of convergence (as is done in~\cite{hs-cardy-embedding}).  
See~\cite{ghs-mating-survey} for a survey of mating-of-trees theory and its applications.

\section{Open problems} \label{sec-open-problems}
 
Here, we discuss some of the most important open problems in the theory of Liouville quantum gravity. 
Additional open problems can be found, e.g., in~\cite[Section 7]{gm-uniqueness}.
Our first problem was alluded to in Section~\ref{sec-lqg-def}. 

\begin{prob} \label{prob-dim}
What is the Hausdorff dimension $d_\gamma$ of a $\gamma$-LQG surface, viewed as a metric space, for $\gamma \in (0,2)\setminus \{\sqrt{8/3}\}$?
\end{prob}

Recall that the Hausdorff dimension of $\sqrt{8/3}$-LQG is 4. 
The value of $d_\gamma$ for general $\gamma\in(0,2)$ is not known even at a physics level of rigor. The best-known physics prediction for $d_\gamma$, due to Watabiki~\cite{watabiki-lqg}, was proven to be incorrect, at least for small values of $\gamma$, by Ding and Goswami~\cite{ding-goswami-watabiki}. 
However, it is known that $\gamma\mapsto d_\gamma$ is strictly increasing~\cite{dg-lqg-dim} and there are reasonably sharp upper and lower bounds for $d_\gamma$~\cite{dg-lqg-dim,gp-lfpp-bounds,ang-discrete-lfpp}. For example, one has $3.550408 \leq d_{\sqrt 2} \leq 3.63299$. 
See also~\cite{bb-lqg-dim} for the most up-to-date simulations of $d_\gamma$. 
Many quantities associated with random planar maps and LQG can be expressed in terms of $d_\gamma$ (see, e.g.,~\cite{dg-lqg-dim,ghs-map-dist,lqg-metric-estimates,gp-kpz,gm-spec-dim,gh-displacement,gp-dla}), so computing $d_\gamma$ would yield many additional results. 

\begin{prob} \label{prob-conv}
Show that weighted random planar map models of the type discussed in Section~\ref{sec-rpm} converge in distribution to $\gamma$-LQG surfaces with $\gamma \in (0,2)\setminus\{\sqrt{8/3}\}$ in both the Gromov-Hausdorff sense and under suitable embeddings into $\BB C$.
\end{prob}

As noted in Section~\ref{sec-conv}, both types of convergence have already been established for certain types of uniform planar maps toward $\sqrt{8/3}$-LQG.

It is possible to construct regularized random measures associated with log-correlated Gaussian fields on $n$-dimensional manifolds for any $n\in\BB N$; see~\cite{rhodes-vargas-review,berestycki-gmt-elementary}.
However, the associated metric has only been constructed in dimension 2.

\begin{prob} \label{prop-higher-dim}
Is it possible to construct random metrics associated with log-correlated Gaussian fields on $\BB R^n$, or on $n$-dimensional manifolds, for $n\geq 3$?
\end{prob}

It is not clear whether random measures and metrics related to log-correlated Gaussian fields in dimension greater than 2 should have any connection to higher-dimensional analogs of random planar maps. 
We remark, however, that a potential analog of the Brownian map in higher dimensions has recently been proposed in~\cite{ml-iterated-folding}. 

\section{Additional expository references}

We mention a few more expository references in addition to the papers cited above.
See~\cite{berestycki-lqg-notes} for introductory lecture notes on the GFF and LQG which go into substantially more detail than this article. 
See~\cite{legall-sphere-survey} for a survey on the geometry of random planar maps and the Brownian map and~\cite{curien-peeling-notes} for lecture notes on random planar maps, emphasizing the applications of the spatial Markov property. 
See~\cite{vargas-dozz-notes} for lecture notes on the conformal field theory (path integral) approach to LQG. 
Although not emphasized in this article, Liouville quantum gravity is closely related to \emph{Schramm-Loewner evolution} (SLE), a family of random fractal curves introduced by Schramm~\cite{schramm0} (see, e.g.,~\cite{shef-zipper,wedges} for relationships between SLE and LQG). 
For an introduction to SLE see the lecture notes~\cite{werner-notes,bn-sle-notes} and the textbook~\cite{lawler-book}.

\medskip
\noindent\textbf{Acknowledgements.} We thank two anonymous referees and also Nina Holden, Jason Miller, and Scott Sheffield for helpful comments. 
We thank J\'er\'emie Bettinelli  and Jason Miller for allowing us to use their beautiful simulations in this article. 
The author was supported by a Clay Research Fellowship and a Junior Research Fellowship at Trinity College, Cambridge. 
\bigskip

\bibliography{cibiblong,cibib}
\bibliographystyle{hmralphaabbrv}


\end{document}